\documentclass[12pt]{article}
\usepackage{graphics}
\usepackage{amsmath}
\usepackage{amssymb}
\usepackage{graphicx}
\usepackage{makeidx}
\usepackage{mathrsfs}
\usepackage{amsfonts}
\usepackage[mathscr]{eucal}
\usepackage{amsmath}
\usepackage{lineno}

\usepackage{graphics}
\DeclareGraphicsExtensions{.eps,.pdf,.jpg,.png}

\def\e{{\rm e}}
\def\e{\hbox{e}}

\def\ds{\displaystyle}

\def\RR{\vbox {\hbox to 8.9pt {I\hskip-2.1pt R\hfil}}}

\def\CC{{\rm C\hskip-4.8pt \vrule height 6pt width 12000sp\hskip 5pt}}

\def\q{\quad}  
\def\cen{\centerline}
\def \rec#1{{\frac{1}{#1}}}
\def\pni{\par\noindent}
\def\vsh{\smallskip}

\def\vsp{\vsh\pni} 

\def\be{\begin{equation}}
\def\ee{\end{equation}}
\DeclareGraphicsExtensions{.eps,.pdf,.jpg,.png}
\def\L{{\mathcal L}} 
\begin{document}

\cen{{\bf FRACALMO PRE-PRINT: \ http://www.fracalmo.org}}

\vsh
\hrule
\vskip 0.25truecm
\font\title=cmbx12 scaled\magstep2
\font\bfs=cmbx12 scaled\magstep1
\font\little=cmr10
\begin{center}
{\title Laplace-Laplace analysis of}
\\[0.25truecm]
{\title the fractional Poisson process}
 \\  [0.25truecm]
 Rudolf GORENFLO $^{(1)}$   and
 Francesco MAINARDI$^{(2)}$
\\[0.25truecm]
$\null^{(1)}$ {\little Department of Mathematics and Informatics,
 Free University of Berlin, }
 \\ {\little Arnimallee  3, D-14195 Berlin, Germany}
 \\{\little E-mail: gorenflo@mi-fu.berlin.de
\\  [0.25truecm]
$\null^{(2)}$ {\little Department of Physics and Astronomy, University of Bologna, and INFN,}
\\
{\little Via Irnerio 46, I-40126 Bologna, Italy}
\\{\little E-mail: francesco.mainardi@unibo.it \  francesco.mainardi@bo.infn.it}
}
\\[0.25truecm]
{ Revised  Version of the paper  
by R. Gorenflo and F. Mainardi,}
\\ {Laplace-Laplace analysis of the fractional Poisson process,}
\\ {in
S. Rogosin   (Editor), {\it Analytical Methods of Analysis and Differential Equations},}
\\
{Belarusian State University, Minsk, 2012, pp. 43--58.}
\\ {[Kilbas Memorial Volume, AMADE 2011]}
\end{center}
\begin{abstract}
\noindent
We generate the fractional Poisson process by subordinating the
standard  Poisson process to the inverse stable subordinator. Our
analysis is based on application of the Laplace  transform with respect to both 
arguments of the evolving probability densities. 
First we give an outline of basic renewal theory, then of the essentials of the classical Poisson process 
and its fractional generalization  via replacement of the exponential waiting time density
 by one of Mittag-Lefler type. 
 Turning our attention to the probability of the counting number 
 of the fractional Poisson process assuming a given value we find in the transform domain
  a formula analogous to the Cox-Weiss formula in the theory of continuous time random walk. 
  This formula contains  for the jump densities (all increments being positive, in fact equal to 1) 
  the Laplace transform instead of the customary Fourier transform. 
 By manipulating this formula we  arrive, after inversion of the transforms, to a subordination
integral involving the inverse stable subordinator. Stochastic
interpretation of this integral leads to the result that the
fractional Poisson process can be obtained from the classical Poisson
process via time change to the inverse stable subordinator.
\end{abstract}
\vskip 1truecm
\centerline{\bf CONTENTS}
\begin{itemize}
\item{Section 1: Introduction}
\item{Section 2: Preliminaries}
\item{Section 3: Elements of renewal theory and CTRW} 
\item{Section 4:  The Poisson process and its fractional generalization}
\item{Section  5:   Subordination}
\item{Section  6:   Conclusions}
\item{References}
\end{itemize}
\section{Introduction}\label{sec:1}
\setcounter{section}{1}
The purpose of this paper is to revisit the fractional Poisson process and present 
its basic theory by treating it as a renewal process with probability of waiting 
time exceeding duration $t$ being given by a Mittag-Leffler type function as 
$E_\beta(-t^\beta)$ with $0<\beta \le 1$. 
Thereby we treat the general renewal process formally as a continuous time random walk (CTRW) 
(stressing this concept) with the counting number playing the role of position in space. 
We then use the known techniques of analyzing the general renewal process and 
its specialization to the Mittag-Lefler waiting time probability distribution,  
working, however, in the transform domain, with the Laplace transform not only 
with respect to time but also with respect to space. 
This is our motivation for calling our method "Laplace-Laplace analysis". 
\vsp
The structure of our paper is as follows.
Section 2 is devoted to notations and terminology, 
in particular to the essential properties of the special functions needed. 
In Section 3 we discuss the elements of renewal theory and the CTRW concept, 
then in Section 4 the Poisson process and its fractional generalization. 
In Section 5 we consider the aspects of subordination, 
thereby also touching the method of parametric subordination for which we cite 
our papers \cite{Gorenflo-Mainardi_EPJ-ST11,Gorenflo-Mainardi_METZLER11,GorMaiViv_CSF07}.
Finally, Section 6 is devoted to conclusions.
\vsp
There are in the literature many papers on the fractional Poisson process where the authors
have outlined a number of aspects and definitions of this process, see e.g.
 Repin and Saichev  
\cite{Repin-Saichev_RQR00}, 
Jumarie  \cite{Jumarie_CSF01}, 
Wang et al.   \cite{Wang-Wen_CSF03,Wang-Wen-Zhang_CSF06},
Laskin  \cite{Laskin_CNSNS03,Laskin_JMP09},
Mainardi et al. \cite{Mainardi-et-al_VIETNAM04}, 
Uchaikin et al.  \cite{Uchaikin-et-al_IJBC08},  
Beghin and Orsingher  \cite{Beghin-Orsingher_EJP09}, 
Cahoy et al.  \cite{Cahoy-et-al_JSPI10},
Meerschaert et al. \cite{M3_FPP},
Politi et al.  \cite{Politi-Scalas_EPL11},
Kochubei  \cite{Kochubei_IEOT11}.
\vsp
To our knowledge the aspect of subordination of the Poisson process 
to the inverse stable subordinator seems to be dealt with only recently 
by Beghin and Orsingher \cite{Beghin-Orsingher_EJP09}  and 
Meerschaert, Nane and Vellaisamy   \cite{M3_FPP}.
 Beghin and Orsingher call the process so generated by subordination the 
 "first form of the fractional  Poisson process" $\,$ (in fact they consider  other kinds of generalization, too)
  whereas the authors of  \cite{M3_FPP} call  it the "fractal Poisson process"$\,$
   and show that it is a renewal process  with the same waiting time distribution as the fractional Poisson process.
\vsp
Our approach, essentially based on the theory of Laplace transforms,
  is alternative to the approaches of these authors and  turns out to be more direct,
     so we expect that the present analysis will be appreciated by applied scientists
   not so well acquainted with the more modern terminology and
    theory of stochastic processes.
\section{Preliminaries}\label{sec:2}
 \setcounter{section}{2} 
 For the reader's convenience  here we present a
brief introduction to the basic notions necessary for our analysis   
of the fractional Poisson process, including the essential elements of integral transforms, 
fractional calculus and special functions.
The notations in these preliminaries follow our earlier papers concerning related topics, see\  
\cite{Gorenflo_PALA10,GAR_Vietnam03,Gorenflo-Mainardi_CISM97,Gorenflo-Mainardi_BAD-HONNEF08,
Gorenflo-Mainardi_EPJ-ST11,Gorenflo-Mainardi_METZLER11,GorMaiViv_CSF07,
Mainardi-et-al_VIETNAM04,Mainardi-Luchko-Pagnini_FCAA01,Mainardi-Mura-Pagnini_IJDE10}. 
For more details on general aspects the interested reader may consult the  treatises 
by Podlubny  \cite{Podlubny_BOOK99},
Kilbas and Saigo \cite{Kilbas-Saigo_BOOK04},
Kilbas, Srivastava and Trujillo  \cite{Kilbas-et-al_BOOK06}, 
Mathai and Haubold  \cite{Mathai-Haubold_BOOK08},
Mathai, Saxena and Haubold  \cite{Mathai-Saxena-Haubold_BOOK-H-2010}, 
Mainardi  \cite{Mainardi_BOOK10}, 
Diethelm  \cite{Diethelm_BOOK10}.    
  \vsp
 \textbf{Fourier and Laplace transforms}
 \vsp
We generously apply the transforms of Fourier and Laplace to admissible functions or generalized
functions defined on $\RR$ or $\RR^+$, respectively.  
  In the following we  use the symbol $\div$ for the juxtaposition of
a function  with its Fourier or Laplace transform.
A look at  the superscript $\,\widehat{\phantom{}}\,$  for the Fourier transform, 
$\,\widetilde{\phantom{}}\,$ for the Laplace transform reveals their relevant 
juxtaposition.  
We use $x$ argument (associated to $\kappa$) for functions Fourier transformed, and 
$x$ or $t$ argument (associated to $\kappa$ or $s$, respectively) for functions Laplace transformed. 
\\
$$  f(x) \,\div\, \widehat f(\kappa) := \int_{-\infty}^{+\infty} \!\! \e^{i\kappa x}\, f(x)\, dx\,,
\q \hbox {Fourier transform}. $$
 $$ f(x) \,\div\, \widetilde f(\kappa) := \int_{0}^{\infty} \!\! \e^{-\kappa x}\,
 f(x)\, dx\,, \q \hbox {Laplace transform}.$$
 $$  f(t) \,\div\, \widetilde f(s) := \int_{0}^{\infty}\!\! \e^{-s t}\,f(t)\, dx\,,
  \q \hbox {Laplace transform}.$$
  \textbf{Convolutions}
  $$ (u*v)(x) := \int_{-\infty}^{+\infty}\!\! u(x-x')\, v(x')\, d x'\,, \q \hbox{Fourier convolution}.$$
  $$ (u*v)(t) := \int_{0}^{t} \!\! u(t-t')\, v(t')\, d t'\,, \q \hbox{Laplace convolution}.$$
  The meaning of the connective $*$
 will be clear from the context. For convolution powers we have:
 $$ u^{*0}(x) = \delta(x)\,, \; u^{*1}(x)= u(x)\,, \; u^{*(n+1)}(x)= (u^{*n}*u)(x)\,,$$
 $$ u^{*0}(t) = \delta(t)\,, \; u^{*1}(t)= u(t)\,, \; u^{*(n+1)}(t)= (u^{*n}*u)(t)\,,$$
 where $\delta$ denotes the Dirac generalized function.
 \vsp
 \textbf{Fractional integral} 
 \vsp
    The \underline{Riemann-Liouville  fractional integral} of order $\alpha>0$,
 for a sufficiently well-behaved  function $f(t)$ ($t\ge 0$), 
  is defined as  generalization
  of the $n$-fold repeated integral, namely 
  $$J^\alpha_t f(t)= 
  \rec{\Gamma(\alpha  )}\,
  \int_0^t (t-\tau )^{\alpha  -1}\, f(\tau )\,{\mbox{d}}\tau \,, \q \alpha>0\,.$$ 
 The {\it semi-group property} is well-known
$$
J_t^\alpha  \,J_t^\beta = J_t^{\alpha  +\beta} = J_t^\beta  \,J_t^\alpha \,,
   \q \alpha \,,\;\beta  \ge 0\,.$$
We have the following Laplace transform pair 
$$ J_t^\alpha  \;f(t) \div
     \frac{\widetilde f(s)}{s^\alpha }\,,\q \alpha  \ge 0\,,  $$
which is the straightforward generalization
of the corresponding formula for the $n$-fold repeated integral.
  \vsp
 \textbf{Fractional derivatives} 
 \vsp
  The  \underline{Riemann-Liouville  fractional derivative} of order $\alpha>0$,
 $D^\alpha _t $, is
 defined as the {\it left inverse operator} of the corresponding fractional integral $J^\alpha _t$.
  Limiting ourselves to fractional derivatives of order $\alpha \in (0,1)$ we have, 
 for a sufficiently well-behaved  function $f(t)$ ($t\ge 0$),   
$$D^\alpha_t  \,f(t) \!:= D_t^1\, J_t^{1-\alpha}\, f(t) =
  {\ds \rec{\Gamma(1-\alpha )}}
  {\ds \frac{{\mbox{d}}}{{\mbox{d}}t}
  \int_0^t 
    \frac{f(\tau)}{(t-\tau )^{\alpha }}\, {\mbox{d}}\tau}\,, \; 
0  <\alpha  < 1 \,, $$ 
while the corresponding \underline{Caputo derivative} is
 $$
	 \begin{array}{ll}
	  {\ds _*D^\alpha_t  \,f(t)} & \!\! := \! J_t^{1-\alpha}\, D_t^1\,f(t) \! = \!
	   {\ds \rec{\Gamma(1-\alpha )}} 
  {\ds  \int_0^t  \frac{f^{(1)}(\tau)}{(t-\tau )^{\alpha}} \,{\mbox{d}}\tau }\\
  & \!\! = \!\!
  {\ds D^\alpha_t f(t) -  f(0^+) \frac{t^{-\alpha }} {\Gamma(1-\alpha )} =\
 D^\alpha_t  \left[ f(t) - f(0^+) \right]} \,.
  \end{array} 
  $$ 
  Both derivatives yield the ordinary first derivative as $\alpha \to 1^-$ but for $\alpha \to 0^+$
  we have
  $$ D^0_t\,f(t)= f(t)\,, \quad  \,_*D^0_t  \,f(t) = f(t) -f(0^+)\,.$$
 We point out the major utility
of the Caputo fractional derivative  
in treating initial-value problems with Laplace transform. We have
$$  \L[\,_*D^\alpha_t\, f(t); s] 
= s^\alpha \widetilde f(s)- s^{\alpha-1}\, f(0^+)\,, \q 0<\alpha\le 1\,. $$
In contrast the Laplace transform of the Riemann-Liouville fractional derivative 
needs the limit at zero of a fractional integral of the function $f(t)$.
   \vsp
 \textbf{Mittag-Leffler and Wright functions}
 \vsp
 \underline{Mittag-Leffler functions}
 $$
 E_\alpha (z) := \sum_{n=0}^\infty \frac{z^n}{\Gamma (\alpha n+1)}
\,,\q \alpha > 0\,, \q z\in \CC\,.$$
The Mittag-Leffler function is  entire of order $1/\alpha$.
The cases $\alpha= 1,2$ are trivial:
$$ 
\begin{cases}
E_1(\pm z)= \exp (\pm z)\,, \\
 E_2\left(+z^2\right) =  \cosh \, (z)\,, \;
   E_2\left(-z^2\right) =  \cos \, (z)\,.
   \end{cases}
    $$
   In the following we will consider, with $0<\beta\le 1$ and $t\ge 0$,
   the Laplace transform pairs 
  $$ 
  \begin{cases}
\Psi(t) \!=\! E_\beta(-t^\beta) \, \div\, {\ds \widetilde \Psi(s)= \frac{s^{\beta-1}}{1+ s^\beta}}\,,
\\
\phi(t)  \!= \!- {\ds \frac{d}{dt}} E_\beta(-t^\beta)
	\, \div \, {\ds \widetilde \phi(s) \!= \!\frac{1}{1+ s^\beta}}\,.
	\end{cases} $$
	It is worth  noting   
the algebraic decay of $\Psi(t)$ and $\phi (t)$ 
as $t\to \infty$:  
$$
{\ds \Psi(t)\sim    \frac{\sin (\beta  \pi)}{\pi}\,\frac{\Gamma(\beta )}{t^\beta }}\,,
\quad 
 {\ds \phi  (t) 
   \sim   \frac{\sin (\beta  \pi)}{\pi}\,\frac{\Gamma(\beta +1)}{t^{\beta +1}}}\,,
 \quad  t \to +\infty\,.
	$$
Furthermore  $\Psi(t) = E_\beta(-t^\beta)$ is the solution of the fractional  relaxation equation
$$ 
        _*D^\beta_t \Psi(t)  = - \Psi(t)\,, \q t\ge 0\,, \quad \Psi(0) =1\,. $$
We refer to \cite{Gorenflo_PALA10,GAR_Vietnam03,Gorenflo-Mainardi_BAD-HONNEF08} 
for the relevance of Mittag-Leffler functions in
theory of continuous time random walk and space-time fractional diffusion
and in asymptotics of power laws.
\vsp
\underline{Wright functions}
 $$ W_{\lambda ,\mu }(z ) :=
   \sum_{n=0}^{\infty}\frac{z^n}{n!\, \Gamma(\lambda  n + \mu )}\,,
 \q \lambda  >-1\,, \,\q \mu \in \CC\,.$$
 We distinguish the Wright functions of 
the {\it first kind} ($\lambda \ge 0$)
and of the  {\it second kind} ($-1<\lambda<0$).
The Wright function is entire of order
$1/(1+\lambda)$ hence
of exponential type only if $\lambda \ge 0$.
The case $\lambda =0$ is trivial since
$  W_{0, \mu }(z) = { \e^{\, z}/ \Gamma(\mu )}\,.$
\vsp
\underline{$M$-Wright function}
$$ M_\nu (z) :=  W _{-\nu , 1-\nu }(-z) \!=\!
 {\ds \sum_{n=0}^{\infty}
 \frac{(-z)^n }{  n!\, \Gamma[-\nu n + (1-\nu )]} } \,, 
  $$
 where $0<\nu<1$. Special cases are 
 $$  M_{1/2}(z) \!=\!
  \rec{\sqrt{\pi}}\, \exp \left(-{\,z^2/ 4}\right),\;
  M_{1/3}(z)  \!=\! 3^{2/3}  {\rm Ai} \left( {z/ 3^{1/3}}\right)\,.
 $$
 Here $Ai$ denotes the Airy
function.
\vsp 
 \underline{Mittag-Leffler functions as Laplace transforms of $M$-Wright functions}
    $$  M_\nu (t) \,\div\,  E_\nu (-s)\,, \; 0<\nu<1\,.  $$
\underline{Stretched Exponentials as Laplace transforms of $M$-Wright functions}
   $$ \frac{\nu }{  t^{\nu +1}}\,  M_\nu \left( 1/{t^\nu } \right)\,\div\,
    \e^{\ds \,-s^\nu}\,, \;  0<\nu <1\,. $$
$$\frac{1}{  t^{\nu}}\,  M_\nu \left( 1/{t^\nu } \right)\,\div\,
    \frac{\e^{\ds\, -s^\nu}}{s^{1-\nu}}\,, \;  0<\nu <1\,. $$
	Note that $\exp (-s^\nu)$ is the Laplace transform of the extremal (unilateral) stable density
	$L_\nu^{-\nu}(t)$, which vanishes for $t< 0$, whereas $\exp (-s^\nu)/s^{1-\nu}$ is related 
	to the Laplace transform  of the Green function of the
	time-fractional diffusion-wave equation. 
\vsp
\underline{The asymptotic representation of the $M$-Wright function}	
\vsp
	Choosing as a variable $t/\nu $ rather than $t$, the computation of the 
asymptotic representation as $t\to \infty$ by the saddle-point approximation
 yields:
$$  M_\nu (t/\nu ) \sim
   a(\nu )\, {\ds t^{(\nu -1/2)/(1-\nu)}}
   \exp \left[- b(\nu)\,{\ds t^{1/(1-\nu)}}\right]\,,
$$
where $$ a(\nu) = \rec{\sqrt{2\pi\,(1-\nu)}}    >0 \,,  \q
  b(\nu) = \frac{1-\nu }{  \nu }    >0 \,.$$
\vsp 
 \underline{The stochastic  relevance of the $M$-Wright function}	
\vsp
 	For the relevance of the $M$-Wright functions in fractional diffusion and related stochastic processes 
	we refer to  formulas shown in the following,  proved in 
	\cite{Mainardi_BOOK10,Mainardi-Luchko-Pagnini_FCAA01,Mainardi-Mura-Pagnini_IJDE10,%
Mainardi-Pagnini-Gorenflo_FCAA03, Pagnini_FCAA13}. 
Let us also point to our pair of complementary papers
\cite{Gorenflo-Mainardi_EPJ-ST11} and \cite{Gorenflo-Mainardi_METZLER11}
on subordination in fractional diffusion processes and our
   contribution \cite{Gorenflo-Mainardi_ Mathai75} in which we analyze some renewal
        processes related and contrasting to the fractional Poisson process.
	\vsp
From the stochastic point of view,
the $M$-Wright function emerges as a natural generalization
of the Gaussian density as Green function for time-fractional diffusion processes and 
of the pure delta-peak drift density $\delta(x-t)$  as Green function for time-fractional drift processes. 
In fact 
\vsp 
-  for the  \underline{time-fractional diffusion equation} with $0<\beta \le 1$: 
$$\, _*D^{\beta}_t u =   \frac{\partial^2 u}{\partial x^2} \,,\; -\infty<x< +\infty\,, \; t\ge 0\,,
$$
we have, assuming  $u(x,0)=\delta(x)$,
$$ u(x,t)\equiv \mathcal{G}_\beta(x,t) 
 = \frac{1}{2} \rec{ t^{\beta/2}}\, M_{\beta/2}
\left(\frac{|x|}{ t^{\beta/2}}\right) \,; $$ 
- for the  \underline{time-fractional drift equation} with $0<\beta \le 1$: 
$$\,_*D^{\beta}_t u(x,t) = -  \frac{\partial }{\partial x} u(x,t) \,,\; -\infty<x< +\infty\,, \; t\ge 0\,, 
$$
we have, assuming $u(x,0)=\delta(x)$,
$$ u(x,t) \equiv G_\beta^*(x,t) = 
 \left\{
  \begin{array}{ll}
    {\ds t^{-\beta}\, M_\beta\left(\frac{x}{t^\beta}\right)}\,, & x>0\,,\\
  0\,, & x<0 \,.
  \end{array}
  \right.
$$
In the former case  the $M$-Wright function can be extended for $1<\beta\le 2$
to represent the Green function of the corresponding time-fractional diffusion-wave equation,
intermediate between the diffusion and wave equations. 
In the last case the $M$-Wright function plays the role of the {\it inverse stable subordinator},
 as it will be explained later on, see Eq. (37) with $x$ replaced by $t_*$. 

\section{Elements of  Renewal Theory  and  CTRW}\label{sec:3}
 \setcounter{section}{3}
\textbf{Definition} of renewal process:  By a {\it renewal process} we mean an infinite  sequence
$0=t_0<t_1<t_2<\cdots$
  of events separated by i.i.d. (independent and identically distributed) random waiting times
  $T_j=t_j-t_{j-1}$, whose probability density $\phi(t)$  is given as a function
   or generalized function in the sense of Gel'fand and Shilov \cite{Gelfand-Shilov_BOOK64}
   (interpretable as a measure) with support on the positive real axis
   $t\ge 0$,  non-negative: $\phi(t)\ge 0$, and normalized: ${\ds \int _0^\infty \!\!\phi(t)\, dt} =1$,
   but not having a delta peak at the origin $t=0$.
   Note that the instant $t_0=0$ is not counted as an event.
\vsp
We are interested in the {\it counting number process} $x=N(t)$ and the sojourn density
$p(x,t)$ for the counting number
\begin{equation}
\label{eq:(counting number)}
  N(t):=
   \hbox{max} \left\{n | t_n \le t \right\}\,, \; n=0,1,2, \cdots
 \end{equation}
  having the value $x$, furthermore in the expectation
 \be 
 \label{eq:(ren)}
 m(t):= \left< N(t)\right>  = \int _0^\infty \!\! x\,p(x,t)\, dx
 \ee 
which is the mean number of positive events in the interval   $[0,t]$
and is called the {\it renewal function}, see e.g. \cite{Ross_BOOK96}.
Then, we ask for the probabilities
\be  
 \label{eq:(probabilities)}
 p_n(t):= P [N(t)=n]\,, \; n=0,1,2, \cdots
 \ee 
\vsp
\textbf{Definition} of continuous time random walk (CTRW):
A {\it continuous time random walk} is an infinite sequence of i.i.d. spatial positions
$0=x_0, x_1, x_2, \cdots$, separated by random jumps    $X_j=x_j-x_{j-1}$,
whose probability density function $w(x)$
is given as a non-negative function or generalized function (interpretable as a measure) 
with support on the real axis
$-\infty <x < +\infty$  and normalized: ${\ds \int _0^\infty \!\! w(x)\, dx} =1$,
 this random walk being subordinated to a renewal process so that we have a random process
 $x=x(t)$    on the real axis with the property $x(t)=x_n$   for $t_n\le t <t_{n+1}$, $n=0,1,2, \cdots$.
  \vsp
     We are interested in the {\it sojourn probability density} $u(x,t)$
	  of a  particle wandering according to the random process $x=x(t)$
	       being in point $x$ at instant $t$.
\vsp
Let us define the following cumulative probabilities related to the probability density function $\phi(t)$
\be
	\Phi(t)  = \int_0^{t+} \!\! \phi (t')\, dt'\,, \quad 	
	\Psi(t)  = \int_{t+}^\infty \!\! \phi (t')\, dt' = 1-\Phi(t)\,.
\ee
For definiteness, we take $\Phi(t)$ as right-continuous, $\Psi(t)$ as left-continuous.
 When the non-negative random variable  represents
 the lifetime of technical systems, it is common to
 call $\Phi(t):= P \left(T \le  t\right) $  the {\it failure probability}
 and $\Psi(t) := P \left(T > t\right)$
  the {\it survival probability}, because $\Phi(t)$ and $\Psi(t)$ are
the respective probabilities that the system does or does not fail
in $(0, t]$. These terms, however,  are commonly adopted for  any
renewal process.
\vsp
 Now, recalling from Section 2 the definition of convolutions,
  we have for the solution $u(x,t)$  the   {\it Cox-Weiss series}, see \cite{Cox_RENEWAL67,Weiss_BOOK94},
  \be  
 \label{eq:(Cox)}
 u(x,t) = \left(\Psi \,*\, \sum_{n=0}^\infty  \phi^{*n}\, w^{*n}  \right)(x,t)
 \ee 
which intuitively says: Before and at instant $t$  there have occurred no jumps
or exactly 1 jump or exactly 2 jumps or ... and if the last jump has occurred at instant $t'$
 the particle is resting there for a duration $t-t'$. 
 \vsp
 In the Fourier-Laplace domain we have
 \be
 \label{eq:(Psi-phi)}
 \widetilde \Psi(s)= \frac {1-\widetilde \phi(s)}{s}\,,
 \ee
  and
\be
\label{eq:(Weiss_FL)}
\begin{array}{ll}
\widehat{\widetilde{u}} (\kappa,s)
&={\ds \frac{1-\widetilde{\phi}(s)}{s} \,\sum_{n=0}^\infty  \left(\widetilde \phi(s)\,\widehat w(\kappa) \right)^n}\\
&={\ds \frac{1-\widetilde{\phi}(s)}{s} \, \frac{1}{1- \widetilde{\phi}(s)\,\widehat{w}(\kappa)}}
\,.
\end{array}
\ee
This is the famous Montroll-Weiss solution formula for a CTRW, see \cite{Montroll-Weiss_JMP65,Weiss_BOOK94}.   
\vsp
In the special situation of  the jump density having support only on the positive semi-axis   $x\ge 0$
it is convenient to replace in this formula
 the Fourier transform by the Laplace transform to obtain the Laplace-Laplace solution
\be
\label{eq:(Cox_LL)}
\begin{array}{ll}
\widetilde{\widetilde{u}} (\kappa,s)
& ={\ds \frac{1-\widetilde{\phi}(s)}{s}\,\sum_{n=0}^\infty  \left(\widetilde\phi (s)\,\widetilde w(\kappa) \right)^n}
\\
& {\ds =\frac{1-\widetilde{\phi}(s)}{s} \, \frac{1}{1- \widetilde{\phi}(s)\,\widetilde{w}(\kappa)}}
\,.
\end{array}
\ee
Henceforth we assume to have this special situation of
{\it jumps only in positive direction} so that we will work with this Laplace-Laplace formula.
	   \vsp
An essential trick of what follows is that we treat renewal processes as continuous time random walks
with waiting time density  $\phi(t)$ and special jump density $w(x)=\delta(x-1)$
corresponding to the fact that the counting number $N(t)$  increases by 1 at each positive  instant $t_n$
so that $x(t) = n$  for $t_n \le t\le t_{n+1}$..
We then have $\widetilde w(\kappa)= \exp(-\kappa)$  and  get for the counting number process  $N(t)$
  the sojourn density	in the transform domain
 \be
\label{eq:(CW)}
\begin{array}{ll}
\widetilde{\widetilde{p}} (\kappa,s)
&={\ds \frac{1-\widetilde{\phi}(s)}{s} \,\sum_{n=0}^\infty  \left(\widetilde \phi (s)\right)^n\,\e^{-n\kappa}} \\
&={\ds \frac{1- \widetilde{\phi}(s)}{s}\,\frac{1}{1-\widetilde{\phi}(s)\, \e^{-\kappa}}}
\,.
\end{array}
\ee	
From this formula we can find formulas for the renewal function   $m(t)$
and the probabilities $p_n(t)=P\{N(t)=n\}$.
Because  $N(t)$ assumes as values only the non-negative integers,
the sojourn density $p(x,t)$  vanishes if $x$ is not equal to one of these,
but has a delta peak of height $p_n(t)$  for $x=n$ ($n=0,1,2,3,\cdots$).
Hence
\be
 \label{eq:(p_n_ser)}
p(x,t)= \sum_{n=0}^\infty p_n(t)\, \delta(x-n)\,.
 \ee
 Rewriting Eq. (\ref{eq:(CW)}), by inverting with respect to $\kappa$, as
 \be
   \sum_{n=0}^\infty \left(\Psi \,*\, \phi^{*n}\right)(t)\, \delta(x-n)\,,
   \ee
   we identify
    \be
\label{eq:(p_n)}
p_n(t)= \left( \Psi\,*\, \phi^{*n}\right) (t)\,.
\ee
According to the theory of Laplace transform we conclude from Eqs.
(\ref{eq:(ren)}) and  (\ref{eq:(p_n_ser)})
\be
\begin{array}{ll}
m(t)&= {\ds - \frac{\partial}{\partial \kappa} \left.\widetilde p(\kappa,t)\right|_{\kappa=0}}
= {\ds \left.\left(\sum_{n=0}^\infty p_n(t)\,n\, \e^{-n\kappa}\right)\right|_{\kappa=0}}\\
 &= {\ds \sum_{n=0}^\infty n\, p_n(t)}\,.
 \end{array}
\ee
a result naturally expected, and
\be
\label{eq:(renf_L)}
\begin{array}{ll}
\widetilde m(s)
&= {\ds \sum_{n=0}^\infty n\, \widetilde p_n(s)}=
{\ds \widetilde\Psi(s)\, \sum_{n=0}^\infty n \, \left(\widetilde \phi(s)\right)^n} \\
&= {\ds \frac{\widetilde \phi(s)}{s\left(1-\widetilde \phi(s)\right)}}\,,
\end{array}
\ee
thereby using the identity
$$ \sum_{n=0}^\infty nz^n = \frac{z}{(1-z)^2}\,, \quad |z|<1\,.$$
Thus we have found in the  Laplace domain the reciprocal pair of relationships
 \be
 \label{eq:(recip)}
  \widetilde m(s)= \frac{\widetilde \phi(s)}{s(1-\widetilde \phi(s))}\,,
  \quad
  \widetilde \phi(s)= \frac{s\,\widetilde m(s)}{1+ s\,\widetilde m(s))}\,,
  \ee
  saying that the waiting time density and the renewal function mutually determine each other uniquely.
  The first formula of Eq. (\ref{eq:(recip)}) can also be obtained  as the value at  $\kappa=0$
   of the negative derivative for $\kappa=0$    of the last expression in Eq. (\ref{eq:(CW)}).
 From Eq. (\ref{eq:(recip)}) follows the reciprocal pair of relationships in the physical domain
      \be
	\begin{cases}  
	m(t) = \int_0^\infty [1 + m(t-t')]\, \phi(t')\, dt'\,, \;\\ \\
	m^\prime (t) =  \int_0^\infty [1 + m^\prime(t-t')]\, \phi(t')\, dt'\,.
	\end{cases}
	  \ee
The first of these equations usually is called the {\it renewal equation}.

\section{The Poisson process and its fractional  \\ generalization}\label{sec:4}
 \setcounter{section}{4}
 The most popular renewal process is the {\it Poisson process}.
  It is characterized by its {\it mean waiting time} $1/\lambda$   (equivalently by its {\it intensity} $\lambda$),
  which is a given positive number, and by its {\it survival probability} $\Psi(t)= \exp (-\lambda t)$
   for  $t\ge 0$, which corresponds to the  {\it waiting time density}
   $\phi(t)= \lambda \, \exp(-\lambda t)$.
	   With $\lambda=1$   we have what we call the {\it standard Poisson process}.
	   \vsp
	   We generalize the standard Poisson process by replacing the
	   exponential function by  a function of Mittag-Leffler type.
	   With $t\ge 0$ and a parameter $\beta\in (0,1]$
	      we take
	\be
	\label{eq:(17)}
	\left\{
	\begin{array}{ll}
	\Psi(t) & \!=  \Psi _\beta(t) = E_\beta(-t^\beta)\,, \\ 
	\phi(t) & \! =  \phi _\beta(t) = - {\ds \frac{d}{dt}} E_\beta(-t^\beta)
	= \beta t^{\beta-1}E^\prime_\beta(-t^\beta)\,.
	 \end{array}
	 \right.
	\ee
	We call the process so defined the {\it fractional Poisson process}.
	To analyze it we go into the Laplace domain where we have
	\be
	\widetilde \Psi(s)= \frac{s^{\beta-1}}{1+ s^\beta}\,, \quad
	\widetilde \phi(s)= \frac{1}{1+ s^\beta}\,.
	\ee
	 If there is no danger of misunderstanding we will not decorate $\Psi$ and $\phi$  with the index $\beta$.
	The special choice $\beta=1$  gives us the standard Poisson process with
	$\Psi_1(t) = \phi_1(t)= \exp(-t)$.
	\vsp
	Whereas the Poisson process has, as is well known, finite mean waiting time, the
	fractional Poisson process ($0<\beta<1$ )  does not have this property. In fact,
	\be
  \!	\langle T\rangle \!=\! \int_0^\infty\!\!t\,\phi(t)\, dt =
\beta \left.{\ds \frac{s^{\beta-1}}{(1+s^\beta)^2}}\right|_{s=0}
	\!=\! \left\{
	\begin{array} {ll}
	1 \,, & \beta=1\,,\\
	\infty \,, & 0<\beta<1\,.
	\end{array}
	\right.
	\ee
Let us calculate the renewal function  $m(t)$.
Inserting $\widetilde \phi(s)= 1/(1+s^\beta)$ into Eq. (\ref{eq:(CW)}) and taking
$w(x)= \delta(x-1)$ as in Section 3,
 we find for the sojourn density of the counting function $N(t)$  the expressions
 \be
 \label{eq:(3.1)}
\widetilde{\widetilde{p_\beta}} (\kappa,s)
=\frac{s^{\beta-1}}{1+s^\beta -\e^{-\kappa}}
 =  \frac{s^{\beta-1}}{1+s^\beta} \,\sum_{n=0}^\infty  \frac{\e^{-n\kappa}}{(1+s^\beta)^n}
\,,
\ee
and
\be
\label{eq:(3.2)}
\widetilde{p_\beta} (\kappa,t) = E_\beta\left(-(1-\e^{-\kappa})t^\beta\right)\,,
\ee
and then
\be
m(t)= - \frac{\partial}{\partial \kappa} \left.\widetilde {p_\beta} (\kappa,t)\right|_{\kappa=0}
= \left. \e^{-\kappa} t^\beta E^\prime_\beta\left( -(1-\e^{-\kappa})t^\beta\right)\right|_{\kappa=0}\,.
\ee
Using $E^\prime_\beta(0)= 1/\Gamma(1+\beta)$
 now yields
 \be
 \label{eq:(3.3)}
 m(t) = \left\{
 \begin{array}{ll}
  t\,, & \beta=1\,, \\
 {\ds \frac{t^\beta} {\Gamma(1+\beta)}}\,, & 0<\beta<1\,.
 \end{array}
 \right.	
 \ee
 This result can also be obtained by plugging $\widetilde\phi(s)= 1/(1+s^\beta)$ into
 the first equation in (\ref{eq:(recip)}) 
 which yields $\widetilde m(s)=1/{s^{\beta+1}}$ and then by Laplace inversion
 Eq. (\ref{eq:(3.3)}).
    \vsp
Using general  Taylor expansion
\be
E_\beta(z) = \sum_{n=0}^\infty \frac{E_\beta^{(n)}}{n!} (z-b)^n\,,
\ee
in Eq. (\ref{eq:(3.2)}) with $b=-t^\beta$  we get
 \be
 \label{eq:(3.3b)}
 \begin{array}{ll}
 \widetilde {p_\beta} (\kappa,t) &=
 {\ds \sum_{n=0}^\infty \frac{t^{n\beta}}{n!}\, E^{(n)}_\beta (-t^\beta)\, \e^{-n\kappa}}\,, \; \\
 p_\beta (x,t) &=
 {\ds \sum_{n=0}^\infty \frac{t^{n\beta}}{n!}\, E^{(n)}_\beta (-t^\beta)\, \delta(x-n)}\,,
 \end{array}
 \ee
and, by comparison with Eq. (\ref{eq:(p_n_ser)}) of Section 3, the probabilities
         \be
		 p_n(t)= P\{N(t)=n\} = \frac{t^{n\beta}}{n!}\, E^{(n)}_\beta (-t^\beta)\,.
		 \ee
Observing from Eq. (\ref{eq:(3.1)})
\be
\widetilde{\widetilde{p_\beta}} (\kappa,s)
=\frac{s^{\beta-1}}{1+s^\beta -\e^{-\kappa}}
= 	  \frac{s^{\beta-1}}{1+s^\beta} \,\sum_{n=0}^\infty  \frac{\e^{-n\kappa}}{(1+s^\beta)^n}\,,
\ee
and inverting with respect to $\kappa$,
\be
\widetilde{p_\beta} (x,s)=
\frac{s^{\beta-1}}{1+s^\beta} \,\sum_{n=0}^\infty  \frac{\delta(x-n)}{(1+s^\beta)^n}\,,
\ee
we finally identify
\be
\label{eq:(3.4)}
\widetilde p_n(s)= \frac{s^{\beta-1}}{(1+s^\beta)^{n+1}} \,\div\,
\frac{t^{n\beta}}{n!}\, E_\beta^{(n)}(-t^\beta) = p_n(t)\,.
\ee
En passant we have proved an often cited special case of an inversion formula 
by Podlubny \cite{Podlubny_BOOK99}.
For an alternative  representation of $p_n(t)$ as an integral see
        Eq. (39) in Section 5.
 \vsp
 For the  Poisson process with intensity $\lambda >0$ we have a well-known infinite system of 
 ordinary differential equations (for $t\ge 0$), see e.g. Khintchine  \cite{Khintchine_QUEUING60},
 \be
 p_0(t)=\e^{-\lambda t}\,,\quad \frac{d}{dt}p_n(t)= \lambda \left(p_{n-1}(t)- p_n(t)\right)\,,
 \ee
   with initial conditions  $p_n(0)=0$, $n=1,2,\dots$
      which sometimes even is used to define the Poisson process.
	  We have an analogous system of fractional differential equations for the fractional Poisson process.
In fact, from Eq. (\ref{eq:(3.4)}) we have
\be
(1+s^\beta)\, \widetilde p_n(s) =
\frac{s^{\beta-1}}{(1+s^\beta)^n}= \widetilde p_{n-1} (s)\,.
\ee
Hence
\be
s^\beta \, \widetilde p_n(s)= \widetilde p_{n-1}(s) - \widetilde p_{n}(s)\,,
\ee
so in the time domain
\be
\label{eq:(33)}
 p_0(t)=E_\beta(-t^\beta)\,,\quad \,_*D_t^\beta p_n(t)= p_{n-1}(t)- p_n(t)\,,
 \ee
with initial conditions  $p_n(0)=\delta_{n\,0}$, $n=0,1,2,\dots$, where
$\,_*D_t^\beta$ denotes the time-fractional derivative of Caputo type of order $\beta$.
It is also possible to introduce and define the fractional Poisson process 
by this difference-differential system.
\vsp 
 Let us note that by solving the system  (\ref{eq:(33)}), Beghin and Orsingher in \cite{Beghin-Orsingher_EJP09}
 introduce what they call the "first form of the fractional  Poisson process"$\,$, and in \cite{M3_FPP}
Meerschaert, Nane and Vellaisamy  show that this process is a renewal process with 
Mittag-Leffler waiting time density as in (\ref{eq:(17)}), hence is identical with 
the "fractional Poisson process".

\section{Subordination}\label{sec:5}
 \setcounter{section}{5}

In order to introduce  the subordination framework in a given stochastic process $x=x(t)$
we must introduce in addition to  the natural time  $t\ge 0$ another time line $t_*\ge 0$,
that we call {\it operational time}. These two basic time lines are inter-related by two 
increasing  processes,  $t=t(t_*)$ and $t_*=t_*(t)$, inverse to each other.  
\vsp
We circumvent the analytical delicacies involved in inverting increasing (but not necessarily  
strictly monotonically increasing) functions with (possibly even in a finite interval infinitely many) 
jumps and intervals of constancy by considering for such functions jumps graphically 
represented by vertical segments and intervals of constancy represented by horizontal segments, 
in the corresponding Cartesian systems of coordinates. 
By inversion (reflection on the diagonal line of the first quadrant) jumps (vertical segments) 
become intervals of constancy (horizontal segments), and vice versa. 
For specialists in stochastic processes let us just  say  that each of the two processes 
(or "trajectories") can (if wanted) be   represented by an equivalent process with c\`adl\`ag structure.
\vsp    
Let us now denote by
$x=x(t)$ the described $\beta$-fractional Poisson process $x=N_\beta(t)$
(happening in $t\ge 0$, running along $x\ge 0$)
and by
$y=y(t_*)$ the standard Poisson process $y=N_1(t_*)$
(happening in  $t_*\ge 0$, running along $y\ge 0$).
Let us denote by 
$t=t(t_*)$ the  totally positive-skewed $\beta$-stable process 
(happening in operational time $t_*$, running along natural time $t$).
Its inverse process 
$t_*=t_*(t)$  (happening in natural time $t$, running along operational  time $t_*$)
is the inverse stable subordinator of order $\beta$.
These two processes $t=t(t_*)$ and $t_*=t_*(t)$ are described in our recent papers
\cite{Gorenflo-Mainardi_EPJ-ST11,Gorenflo-Mainardi_METZLER11},
 where they are referred to as the {\it leading process} and {\it directing process}, respectively.
\vsp
Thus, from now on, we will consider the Poisson process on the time line $t_*$,
that we call {\it operational time}, and we  denote it by $y=y(t_*)$.
Its density in $y$, evolving in operational time $t_*$,  is
$p_1(y,t_*)$.
We will now state a  theorem on time change and derive its underlying  subordination integral.
\vsp
\textbf{Theorem on Subordination and Time Change}
\be
\label{eq:(10)}
x(t)= y(t_*(t))\,.
\ee
The {\bf proof} will be carried via the Laplace-Laplace method.
\vsp
By Eq.(\ref{eq:(3.1)}) with $\beta=1$ 
we have for the standard Poisson process (happening in operational time $t_*\ge 0$
and running on $y\ge 0$)
$$ \widetilde{\widetilde {p_1}}(\kappa,s_*)= \frac{1}{s_* +1 -\exp(-\kappa)}\,.$$
Then the first Laplace inversion $s_* \to t_*$ yields by by Eq.(\ref{eq:(3.2)}), or directly,
$${\widetilde {p_1}}(\kappa,t_*) = \exp(-t_*(1-\e^{-\kappa}))\,.$$
and the second Laplace inversion $\kappa\to y$ yields by Eq.(\ref{eq:(3.3b)}), or directly,
\be
\label{eq:(NEW)}
p_1(y,t_*)  ={\ds \sum_{n=0}^\infty \frac{t_*^{n}}{n!}\, \exp ( -t_*) \, \delta(y-n)}\,.
\ee
 Using
$$    \int_0^\infty
\exp(-at_*)\, dt_* = \frac{1}{a} \quad \hbox{for} \quad
a=s^\beta +1 -\e^{-\kappa}\,,$$
we write by Eq.(\ref{eq:(3.1)}) in the form
  $$\widetilde{\widetilde {p_\beta}}(\kappa,s)=
  s^{\beta-1} \int_0^\infty
  \exp [-t_*(s^\beta+1-\e^{-\kappa})] \,  dt_*$$
  $$ = \int_0^\infty
  \exp\left( -t_*(1-\e^{-\kappa})\right)
  \{ s^{\beta-1} \exp (-t_*s^\beta)\} \,dt_*\,.$$
Then,  by Laplace-Laplace inversion, we obtain the \textbf{subordination integral}
\be
\label{eq:(11)}
 p_\beta (x,t) = \int_0^\infty p_1(x,t_*)\, q_\beta(t_*,t)\, dt_*\,,
\ee
with the density  $q_\beta(t_*,t)$  of the inverse stable subordinator
(density in operational time $t_*$ evolving in natural time $t$)
that we have used in \cite{GorMaiViv_CSF07}
and  later  in  \cite{Gorenflo-Mainardi_EPJ-ST11,Gorenflo-Mainardi_METZLER11}.
In particular, from our  most recent paper \cite{Gorenflo-Mainardi_METZLER11}, see Eqs. (75) and (81) in it,
we have
\be
\label{eq:(36)}
   q_\beta(t_*,t) = \, J_t^{1-\beta}\, r_\beta(t, t_*) = t^{-\beta}\, M_\beta(t_*/t^\beta)
\, \div\, s^{\beta-1} \, \exp (-t_*s^\beta)    \,.
\ee
Here
$\, J_t^{1-\beta}$ denotes the Riemann-Liouville time-fractional integral of order $1-\beta$,
  $r_\beta(t, t_*)$ is the probability density (in $t$, evolving in operational time $t_*$) of
the extremely positively skewed stable density of index $\beta$, $L_\beta^{-\beta}$, namely
\be
\label{eq:(37)}
 r_\beta(t,t_*) = t_*^{-1/\beta} L_\beta^{-\beta} (t/t_*^{1/\beta})
\div \e^{-t_*s^\beta}\,,
\ee
and  $M_\beta$ is the so-called $M$-Wright function of order $\beta$ defined 
in Section 2.
Eqs. (\ref{eq:(36)}), (\ref{eq:(37)})
appear in various disguises in publications
 of other authors for which we have not done a systematic search, see e.g.  \cite{M3_PRE02}.
\vsp
\textbf{Interpretation}: Equation (\ref{eq:(11)}) says:
The fractional Poisson process $x(t)$  can be obtained from the standard Poisson process
by time change via the inverse stable subordinator $t_*(t)$ according to Eq.(\ref{eq:(10)}), namely
 $ x(t)= y(t_*(t))$.
\vsp \vsp
 \textbf{Corollary 1}: For the fractional Poisson probabilities we have, alternatively to (26),
 the integral representation (the subordination integral)
 \be 
 \label{eq:(NEW2)}
p_n(t)\!= \! \frac{1}{n!} \, \int_0^\infty \!\!  t_*^n \exp(-t_*) \, q_\beta(t_*,t)\, dt_* \!=\!
   \frac{1}{n!} \, \int_0^\infty\!\!  t_*^n\, t^{-\beta} \exp(-t_*) \, M_\beta(t_* t^{-\beta})\, dt_*.
\ee 
\textbf{Proof}: Set $x=n$ in (35), use (26) with $E_1(-t_*) =\exp(-t_*)$, (35) and (36). 
\vsp
As an exercise the reader may verify  (39) directly from  (26) by use of the fact that
      $E_\beta(-s)$ is the Laplace transform of $M_\beta(t)$   (see equation (F.30) in [20]).
 \vsp \vsp
\textbf{Corollary 2}: With the positive-oriented extremal $\beta$-stable process
$t=t(t_*)$   we have the parametric representation 
\be
\left\{
\begin{array}{ll}
t = t(t_*)\,, \\
x = y(t_*)\,,\\
\end{array}
\quad \hbox{for} \quad x=x(t)\,.
\right.
\ee
\vsp
\textbf{Remark}
According to Meerschaert, Nane and Vellaisamy  \cite{M3_FPP}
\be
z_\beta (t) = z_1 (y(t_*(t)))\,,
\ee
where $z_\beta (t)$ is a CTRW with $\beta$-Mittag-Leffler
waiting times and an arbitrary jump density, and $z_1$ is a CTRW with the exponential
waiting times of the standard Poisson process having the same arbitrary jump distribution.
The influence of the index $\beta$ there is transferred to the process $t_*(t)$.
\vsp
IN WORDS: Define a random walk as the sum of $N$ spatial random variables,
$N$ being the counting number of a renewal process.
Then, the subordination  of such random walk to the fractional Poisson process
is equivalent to the replacement of the fractional Poisson process by the standard Poisson process
evolving on the inverse stable subordinator of order $\beta$.
\section{Conclusions}
 By time change via the inverse stable subordinator the
standard Poisson process is transformed to the fractional Poisson process.
Treating the counting number of a renewal process as a "spatial" $\,$ variable and not 
shying away from extensive use of delta-functions and 
their transforms we can use the formalism of the the common theory of Continuous Time Random Walk. 
In this theory usually in space the Fourier transform is used, but as in our situation all jumps are 
of size 1 our walk proceeds only in positive direction, 
hence not only in time but also in space we can work with the transform of Laplace instead 
of that of Fourier. 
\vsp
The analysis carried out in this essay is based on the Laplace transform of the waiting time density 
as well as of the jump-width density. 
So, in the Montroll-Weiss solution formula written in Laplace-Laplace instead of Laplace-Fourier form, 
we apply the usual trick of representing the reciprocal of the denominator as an improper integral 
(thereby introducing an operational time variable) 
and can (by inverting the transforms) separate variables of time and space.
Finally, we  arrive at an integral representation that allows interpretation as a time change realized 
(as in the theory of time-fractional diffusion) by the "inverse stable subordinator".

\vspace{3mm}


\end{document}